\documentclass[11pt,a4paper]{article}
\vspace{7cm} \textwidth 14cm \textheight 21.6cm
\usepackage{amsmath}
\usepackage{amsmath}
\usepackage{amsfonts}
\usepackage{latexsym}
\usepackage{latexsym,amsthm,amscd}
\newtheorem{thm}{Theorem}

\newtheorem{lem}{Lemma}
\newtheorem{prop}{Proposition}
\newtheorem{defn}{Definition}

\numberwithin{equation}{section}
\newcounter{alphthm}
\newtheorem{prp}[alphthm]{Theorem}
\newcommand{\be}{\begin{equation}}
\newcommand{\ee}{\end{equation}}
\newcommand{\ben}{\begin{enumerate}}
\newcommand{\een}{\end{enumerate}}
\newcommand{\beq}{\begin{eqnarray}}
\newcommand{\eeq}{\end{eqnarray}}
\newcommand{\beqn}{\begin{eqnarray*}}
\newcommand{\eeqn}{\end{eqnarray*}}

\newcommand{\bl}{\begin{lem}}
\newcommand{\el}{\end{lem}}
\newcommand{\bp}{\begin{prop}}
\newcommand{\ep}{\end{prop}}
\newcommand{\bd}{\begin{Def}}
\newcommand{\ed}{\end{Def}}
\newcommand{\bt}{\begin{thm}}
\newcommand{\et}{\end{thm}}
\begin{document}
\title{COMPLETE LIFT CONFORMAL VECTOR FIELDS ON FINSLER MANIFOLDS
\thanks{The author should express his gratitude to prof. Akbar-Zadeh for
his valuable suggestions.}}
\author{B. BIDABAD}

\date{ }
\maketitle

\begin{abstract}

In Finsler geometry the complete lift vector fields have
distinguished geometric significance. For example a vector field on
a Finsler manifold is said to be conformal if its complete lift is
conformal in usual sense. In this work  we define a new Riemannian
or Pseudo-Riemannian metric on TM derived from a Finsler metric on
the base manifold $M$. This metric is in some senses more general
than the other lift metrics defined previously on TM and then  we
study the complete lift vector fields on $TM$. More precisely we
prove; \emph{Let $(M,g)$ be a Finsler manifold, $TM$ its tangent
bundle and \~g a Riemannian metric on TM derived from g. Then every
complete lift conformal vector field is homothetic}.
\end{abstract}
 \vspace{-5in}
 {\it \small {\it \small Balkan Journal of Geometry and Its Applications, Vol.11, No.1, 2006.}}
\vspace{5in}\\
 {\small{\textbf{Mathematics Subject
Classification}2000: 53C20. }}
 \section*{Introduction.}
 Let $(M,g)$ be an n-dimensional Riemannian manifold and
 $\phi$ a transformation on $M$. Then
 $\phi$ is called a {\it conformal} transformation, if it  preserves the angles.
 Let $V$ be a vector field on $M$ and  $\{\varphi_{t}\}$ be the local
 one-parameter group of local transformations on $M$
 generated by $V$. Then $V$ is called a {\it
  conformal vector field}
 on $M$ if each $\varphi_{t}$ is a
 local  conformal transformation of $M$.
 It is well known that $V$ is a {\it conformal vector
 field} on $M$ if and only if there is a scalar function $\rho$ on $M$ such that
 $\pounds_{_{V}}g=$2${\rho}g$ where $\pounds_{_{V}}$ denotes  Lie derivation
 with respect to the vector field $V$. Specially $V$ is called
  {\it  homothetic} if $\rho$ is constant
  and it is called an {\it isometry } or {\it Killing vector field}
  when $\rho$ vanishes.\\

  There are some lift metrics on $TM=\cup_{_{x\in M}} T_xM$ as follows:
    {\it complete} lift metric or  $g_ 2$,
   {\it diagonal} lift metric or   $g_ 1+g_ 3$,
   lift metric $g_ 2+g_ 3$
  and   lift metric $g_ 1+g_ 2$, where $g_{1}:=g_{ij}dx^{i}\otimes dx^{j}$,
 $g_{2}:=2g_{ij}dx^{i}\otimes\delta{y^{j}}$ and
 $g_{3}:=g_{ij}\delta{y^{i}}\otimes\delta{y^{j}}$ are all bilinear differential forms defined
 globally on $TM$.\\

 In the study of Finsler geometry the complete lift vector fields have a great significance.
 More precisely  let  $V$ be a vector field on the Finsler manifold
  $(M,g(x,y))$ and $  X^c$
  be the complete lift of $V$. Then $V$ is called a  \emph{conformal vector field  of
  Finsler manifold }$(M,g)$ if there is a scalar function\footnote{By a simple calculation
and vertical partial derivative using commutative property of Lie
derivative one can show that $\Omega $ is a function of $x$ alone
[Ak].}
 $\Omega$ on $TM$  which satisfies $\pounds_{_{ X^c}} g=2\Omega g.$\\
For the complete lift vector fields the following results are well
known:
\newline
\begin{prp}
{[YI]:}Let $( M, g)$ be a Riemannian manifold, $X$  a vector field
 on $M$ and  $X^{C}$  complete lifts of $X$ to $TM$. If we consider  $TM$ with  metric
 $g_ 2$ then $X^{C}$ is a conformal vector field on $TM$  if and only if
$X$ is homothetic on $M$.
\end{prp}

\begin{prp}
{[YK]:}Let $( M, g)$ be a Riemannian manifold. If we consider $TM$
with metric $g_ 1+g_ 3$ then $X^{C}$ is a conformal vector field on
$TM$ if and only if $X$ is homothetic.
\end{prp}

 In a recent work we introduced  a new Riemannian and
 pseudo-Riemannian lift
  metrics on $TM$, $\widetilde{g}=ag_1+bg_2+cg_3$ where $a$, $b$ and $c$ are
certain constant real numbers. That is a combination of diagonal
lift, and complete lift metrics, which is in some senses more
general than those who are used previously.  We have replaced the
cited lift metrics in Theorems A and B by $\widetilde{g}$.
  More precisely, we have proved Theorem C in [BH] as follows.

 \begin{prp}Let $M$ be an  n-dimensional
Riemannian manifold and let $TM$ be its tangent bundle with  metric
$\widetilde{g}$. Then every complete lift conformal vector field
 on $TM$ is homothetic.
 \end{prp}

In the present work we replace the Riemannian metric on $M$ by a
Finsler metric endowed with a Cartan connection and prove the  following theorem. \\
 \newline  {\bf Theorem 1}\textit{: Let (M,g) be a $C^{\infty}$ Finsler manifold with
Cartan connection, $TM$ its tangent bundle and \~g the Riemannian
(or Pseudo-Riemannian) metric on TM derived from g. Then every
complete lift conformal vector field on $TM$ is homothetic.}\\

\section{ Preliminaries.}

  Let $M$ be a real n-dimensional  manifold of class $C^ \infty$. We
denote by  $TM\rightarrow M$ the  bundle
  of tangent vectors and by $ \pi:TM_{0}\rightarrow M$ the fiber bundle of
non-zero tangent vectors.
  A {\it{Finsler structure}} on $M$ is a function
$F:TM \rightarrow [0,\infty )$, with the following properties: (I)
$F$ is differentiable ($C^ \infty$) on $TM_{0}$; (II) $F$ is
positively homogeneous of degree one in $y$, i.e.
 $F(x,\lambda y)=\lambda F(x,y),  \forall\lambda>0$, where we denote
 an element of $TM$ by $(x,y)$.
(III) The Hessian matrix of $F^{2}$ is positive definite on
$TM_{0}$; $(g_{ij}):=\left({1 \over 2} \left[
\frac{\partial^{2}}{\partial y^{i}\partial y^{j}} F^2
\right]\right).$
 A \textit{Finsler
manifold} is a pair of a differentiable manifold $M$ and a Finsler
structure $F$ on $M$. The tensor field $g=(g_{ij})$ is called the
\emph{Fundamental Finsler tensor} or \emph{ Finsler metric tensor}.
Here, we denote a Finsler manifold by $(M,g)$.\\

Let $V_vTM=ker\pi_{*}^v$ be the set of the vectors tangent to the
fiber through $v\in TM_0$. Then a \emph{vertical vector bundle} on
$M$ is defined by $VTM := \bigcup_{_{v\in TM_0}}V_vTM$. A
\textit{non-linear connection} or a \textit{horizontal distribution}
on $TM_0$ is a complementary distribution $HTM$  for $VTM$ on
$TTM_0$. Therefore we have the decomposition
 \be \label{decomp}
TTM_0 =VTM\oplus HTM.
 \ee
$HTM$ is a \emph{ vector bundle} completely determined by the
non-linear differentiable functions
   $N^{j}_{i}(x^{},y^{})$ on $TM$, called coefficients of the
non-linear connection.  Let $HTM$ be a non- linear connection on
$TM$ and $\nabla$ a linear
  connection on $VTM$, then the pair $(HTM,\nabla)$ is called a
  \textit{Finsler connection} on the manifold $M$. \\
Using the local coordinates $(x^{i},y^{i})$ on $TM$ we have the
local field of frames $\{\frac{\partial}{\partial
x_{i}},\frac{\partial}{\partial y_{i}}\}$ on $TTM$.
  It is well known that we can choose a local field of frames $\{\frac{\delta}{\delta x_{i}},\frac{\partial}
  {\partial y_{i}}\}$ adapted to the above decomposition i.e. $\frac{\delta}{\delta x_{i}}\in {\Gamma}(HTM)$
   and $\frac{\partial}{\partial y_{i}}\in {\Gamma}(VTM)$ set of vector fields on $HTM$
  and $VTM$, where
  $ \label{delta}
  \frac{\delta}{\delta x_{i}}:=\frac{\partial}{\partial x_{i}}-N^{j}_{i}\frac{\partial}
  {\partial y_{j}},$ where we use the {\it Einstein summation convention}.

Here, in this paper, all manifolds are supposed to be connected.\\

Let $(M,g(x,y))$ be a Finsler manifold then a Finsler connection is
called a \textit{metric Finsler connection} if $g$ is parallel with
respect to $\nabla$. According to the Miron terminology this means
that $g$ is both horizontally and vertically metric. The
\textit{Cartan connection} is a metric Finsler connection for which
the Deflection , horizontal and vertical torsion tensor fields
vanishes.

 Let $(M,g(x,y))$ be a Finsler manifold with metric
Finsler connection the \textit{Curvature tensors} of $M$ are defined
by
   $$ R(X,Y)Z = \{{[\nabla_{X},\nabla_{Y}]-\nabla_{[X,Y]}}\}Z, $$
where $X , Y , Z \in {\cal X}(TM)$\newline They are called
accordingly to the choice of $X$ and $Y$ in $HTM$ or $VTM$
horizontal or
vertical curvature  tensors of Finsler manifold.\\

 Let $M$ be a Finsler manifold and $\nabla$ a Finsler connection on
$M$, then we have [MA]\\

$R_{k\ ji}^{ \ h}=\delta_{i}F_{k \ \ j}^{ \ h}-\delta_{j}F_{k \ \
i}^{\ h} +F_{k \ \ j}^{\ m}F_{m \ \ i}^{ \ h}-F_{k\ \ i}^{ \ m}F_{m
\ \ j}^{\ h}
+C_{k\ \ m}^{ \ h}R_{j \ \ i}^{ \ m}, $\\

$R^{h}_{ \ \ ij}=\delta_{j}N_{i}^{h}-\delta_{i}N_{j}^{h},$ where we
have put $\partial_{i} = \frac{\partial}{\partial x^{i}} ,
\dot{\partial}_{i}=\frac{\partial}{
\partial y^{i}}$ ,
$\delta_{i}=\partial_{i}-N_{i}^{m}\dot{\partial}_{m}.$\\
If $\nabla$ is a
 Cartan connection then
$N_{i}^{h}=y^{m}F_{m \ i}^{\ h}.$

 \begin{prop}{[Mi]} Let $M$ be an n-dimensional Finsler
space with a Cartan connection, then we have the following equations
\begin{eqnarray*}
&&(1) \ \ F_{i \ \ j}^{\
h}=\frac{1}{2}g^{hm}(\delta_{i}g_{mj}+\delta_{j}g_{im}
-\delta_{m}g_{ij}). \hspace*{9cm}\\
&&(2) \ \ C_{ijk}=\frac{1}{2}\dot{\partial}_{k} g_{ij}\ \ \  where \
\ C_{ijk}=
C_{i \ \ k}^{ \ m} g_{jm}.\\
&&(3) \ \ y^{m}C_{mij}=0.\\
&&(4) \ \ R^{h}_{ \ \ ij}=y^{m}R_{m \ \ ij}^{ \ \ h}.\\
\end{eqnarray*}
\end{prop}
The Cartan horizontal and vertical covariant derivative of a tensor
field of type $(_{2}^1)$ are given locally as follows:\\
\begin{equation}\label{h-der}
\nabla_{j} T_{k\ \ i}^{\ h}:=T_{k\ \ i|j}^{\ h}=\delta_{j}T_{k\ \
i}^{ \ h}+ F_{m\ j}^{ \ h} T_{k\ \ i}^{ \ m}-F_{k\ \ j}^{ \ m} T_{m
\ \ i}^{ \ h} -F_{i\ \ j}^{ \ m} T_{k\ \ m}^{ \ h}.
\end {equation}
\begin{equation*}\label{v-der}
 \nabla_{\bar j} T_{k\ \ i}^{\ h}:=T_{k\ \ i|\overline{j}}^{\ h}=\dot{\partial}_{j}T_{k\
\ i}^{ \ h}+ C_{m\ j}^{ \ h} T_{k\ \ i}^{ \ m}-C_{k\ \ j}^{ \ m}
T_{m \ \ i}^{ \ h} -C_{i\ \ j}^{ \ m} T_{k\ \ m}^{ \ h}.
\end {equation*}
\section{Lift metrics and conformal vector fields.}
 \subsection{Complete lift vector fields and Lie derivative.}
 Let $V=v^i\frac{\partial}{\partial x^i}$ be a vector field on $M$. Then
   $V$ induces an infinitesimal point transformation
 on $M$. This is naturally extended to a point transformation
    of the tangent bundle $TM$ which is called {\it extended point
    transformation}.
    Let $V$ be a vector field on $M$ and $\{\Phi_{t}\}$
  the local one parameter groups of $M$ generated by $V$. Let $\tilde\Phi_{t}$ be
  the extended point transformation of
   $\Phi_{t}$ and $\{\tilde\Phi_{t}\}$ be the local one-parameter groups of $TM$.
  If $ X^c$ is a vector field on $TM$ induced by $\{\tilde\Phi_{t}\}$ it is
  called the \textit{ complete lift} vector field of $V$.\\
    It can be shown that the extended point transformation  is a
    transformation induced by the complete lift vector field of $V$,
    $X^c=v^{i}\delta_i
    +y^{j}\nabla_{j}v^{i}\dot{\partial} _i$
     with respect to the decomposition
    (\ref{decomp}). \\

   Let $M$ be an n-dimensional  manifold, $V$  a
  vector field on $M$ and $\{\phi_{t}\}$  a 1-parameter  local group of
  local transformations of $M$ generated by $V$. Take any tensor
  field $S$ on $M$, and denote by $ {\phi_{t}}^{*}(S)$ the
  pulled back of $S$ by $\phi_{t}$. Then the {\it Lie derivation} of $S$
  with respect to $V$ is a tensor field $\pounds_{_{V}} S $ on $M$
  defined by: \newline
   $$\pounds_{_{V}} S=\frac{\partial}{\partial
   t}{\phi_{t}}^{*}(S)|_{t=0} = \lim _{t\longrightarrow {0}} \frac{ {\phi_{t}}^{*}(S)-(S)}{t},$$ on the domain of $\phi_{t}$.
   The mapping $\pounds_{_{V}}$ which map S to  $\pounds_{_{V}}(S)$ is
   called the Lie derivation with respect to $V$.\newline

In  Finsler geometry the Lie derivative of an arbitrary tensor,
$T^{{\,\,\,\,\,}k}_{ij}$ is given locally by [Yan1]:
$$\pounds_{_V}T^{{\,\,\,\,\,}k}_{i}=v^a \nabla_a T^{{\,\,\,\,\,}k}_{i}+
v^a \nabla_a v^b {\nabla}_{\overline{b}} T^{{\,\,\,\,\,}k}_{i}-
T^{{\,\,\,\,\,}a}_{i}\nabla_av^k + T^{{\,\,\,\,\,}k}_{a}\nabla_iv^a
,$$ or equivalently,

\begin{equation}\label{Lee der}
 \pounds_{_V}T^{\,\,j}_i=v^a\partial_a
T^{\,\,j}_i+y^a\partial_av^b\dot{\partial}_b T^{\,\,j}_i -
T^{\,\,a}_i \partial_av^j+T^{\,\,j}_a\partial_iv^a .
\end {equation}\\
So we have

\begin{equation}\label{Lee der y}
\pounds_{_V}y^i=v^a\partial_ay^i+y^b\partial_bv^j\dot{\partial}_jy^i-y^a\partial_av^i
 =y^b\partial_bv^i-y^a\partial_av^i=0 ,
\end {equation}\\
\begin{equation}\label{Lee der g}
\pounds_{_V}g_{ij}=v^a\partial_a
g_{ij}+y^a\partial_av^b\dot{\partial}_b g_{ij} + g_{aj}\partial_i
v^a+g_{ia}\partial_jv^a.
\end {equation}

 We have also this interchanging formula between  Cartan covariant derivatives and Lie
 derivatives.
 \begin{equation}\label{interchange}
   \nabla_{k}\pounds_{_V} g_{ij}- \pounds_{_V}\nabla_{k}
 g_{ij}=g_{aj}\pounds_{_V}F^{a} _{ik}+ g_{ai}\pounds_{_V}F^{a}
 _{jk}.
 \end{equation}

  \subsection{A lift metric on tangent bundle.}
Let $(M,g)$ be a Finsler manifold. In this section we define a new
Riemannian or Pseudo-Riemannian metric on $TM$ derived from the
Finsler metric. This metric is in some senses more general than the
other lift metrics defined previously on $TM$. By mean of the dual
basis $\{dx^i, \delta y^{i}\}$ analogously to the Riemannian
geometry the tensors;
  $g_{1}:=g_{ij}dx^{i} \otimes dx^{j}$\, \
 $g_{2}:=2g_{ij}dx^{i}\otimes \delta{y^{j}}$\, \,\,
 $g_{3}:=g_{ij}\delta{y^{i}}\otimes \delta{y^{j}}$ \,\,
 are all quadratic differential tensors defined globally on $TM$, see [YI].
 Now let's consider the Finsler metric tensor $g$ with the
 components $g_{ij}(x,y)$ defined on $TM$. The tensor field
$\widetilde{g}=\alpha g_1+\beta g_2+\gamma g_3$ on $TM$, where the
coefficient $\alpha , \beta , \gamma$ are real numbers, has the
components
$$\left(%
\begin{array}{cc}
  \alpha g & \beta g \\
  \beta g & \gamma g \\
\end{array}%
\right)$$ with respect to  the dual basis of $TM$. From the linear
algebra we have $ det \widetilde{g}=(\alpha \gamma -\beta^2)^n det
g^2$. Therefore $\widetilde{g}$  is nonsingular if $\alpha \gamma
-\beta^2\neq0 $ and it is positive definite if $\alpha , \gamma$ are
positive and $\alpha \gamma -\beta^2 > 0 $. Indeed $\widetilde{g}$
define respectively a Pseudo-Riemannian or a Riemannian lift metric
on $T(M)$.
\begin{defn} Let $(M , g)$ be a Finsler manifold. Consider  tensor field $\widetilde{g}=\alpha g_1+\beta g_2+\gamma g_3$ on
$TM$, where the coefficient $\alpha , \beta , \gamma$ are real
numbers. If $\alpha \gamma -\beta^2\neq0 $ then $\widetilde{g}$ is
non-singular and it can be regarded as a  Pseudo-Riemannian metric
on $TM$. If $\alpha$ and $ \gamma$ are positive such that $\alpha
\gamma -\beta^2 > 0 $  then $\widetilde{g}$ is positive definite and
consequently can be regarded as a
Riemannian  metric on $TM$.\\
$\widetilde{g}$ is called the lift metric of $g$ on $TM$.
 \end{defn}

\subsection{Conformal vector fields.}
\renewcommand{\theequation}
{\arabic{section} .\arabic{equation}}
   Let $(TM,G(x,y))$ be a
  Riemannian ( or pseudo-Riemannian) manifold.  A vector field $ \widetilde{X}\in {\cal X}(TM)$ on $TM$
  is called  a
   \it conformal vector field on TM \rm if there is a scalar function $\Omega $ on $TM$
   such that
  \begin{equation}
   \pounds_{_{\widetilde{X}}}G = 2\Omega G.
   \nonumber
  \end{equation}
  If $\Omega$ is constant then the vector field $\widetilde{X}$ is called \it homothetic \rm and
  if $\Omega$ is zero then its called   an \it isometric \rm or a \it Killing \rm vector field \rm. \\
    Now let we consider $(TM,\tilde{g}(x,y))$ with the
  complete
  lift vector field $X^c$ of an arbitrary vector field $V$ on $M$. Then by above
  definition we call ${X}^c$
   a conformal vector field on $TM$ if
   $$ \pounds_{_{{X}^c}}\tilde g = 2\Omega \tilde g.$$
\section{Main results}
\renewcommand{\theequation}
{\arabic{section} .\arabic{equation}}
\setcounter{equation}{0}
Analogous to the Riemannian geometry [Yam], by straight forward
calculation we have  the following lemmas in Finsler geometry.
\begin{lem} : Let $(M,g)$ be a Finsler manifold with Cartan
connection, then we have;\\ (1)  $[X_i , X_j]=
R^h_{\,\,\,\,ij}X_{\overline{h}},$
\newline (2) $[X_i ,
X_{\overline{j}}]=\dot{\partial}_jN^h_{\,\,i}X_{\overline{h}},$
\newline (3)  $[X_{\overline{i}} , X_{\overline{j}}]= 0.$\\
where we put $X_i=\delta_i$ and
$X_{\overline{i}}=\dot{\partial_{i}}$ for simplicity.

\end{lem}
 Let's denote by  $\pounds_{_{X^c}}$ the lie derivative with respect to
the complete lift vector field $X^c$. Then we obtain the following
lemma :\newline
\begin{lem} : Let $(M,g)$ be a Finsler manifold with Cartan
connection, then we have;\\
$(1)
\pounds_{_{X^c}}X_i=-\partial_iv^hX_h-\pounds_{_V}N^h_{\,\,i}X_{\overline{h}},$\newline
$(2)
\pounds_{_{X^c}}X_{\overline{i}}=-\partial_iv^hX_{\overline{h}},$\newline
$(3)   \pounds_{_{X^c}}dx^h=\partial_mv^hdx^m ,$\newline $(4)
\pounds_{_{X^c}}\delta y^h=\pounds_{_V}N^h_{\,\,m}dx^m+\partial_m
v^h\delta y^m.$
\end{lem}

\begin{proof} (1)   $\pounds_{_{X^c}}X_i=[X^c ,
X_i]$\newline
${\,\,\,\,\,\,\,\,\,\,\,\,\,\,\,\,\,\,\,\,\,\,\,\,\,\,\,}=[v^hX_h+y^m
v^h_{\,\,\,\,|m}X_{\overline{h}} {\,\,\,\,}, X_i]$\newline
${\,\,\,\,\,\,\,\,\,\,\,\,\,\,\,\,\,\,\,\,\,\,\,\,\,\,\,}=v^h[X_h ,
X_i]-X_i(v^h)X_h+y^m v^h_{\,\,\,\,|m}[X_{\overline{h}} , X_i]-
X_i(y^m v^h_{\,\,\,\,|m})X_{\overline{h}}$\newline
${\,\,\,\,\,\,\,\,\,\,\,\,\,\,\,\,\,\,\,\,\,\,\,\,\,\,\,}=-\partial_iv^hX_h-\pounds_{_V}
N^h_{\,\,i} X_{\overline{h}}.$\newline Thus we get (1). We can prove
(2) by the same way as the proof of (1). Next we prove (3). Since
$(dx^h , \delta y^h)$ is the dual basis of $(X_h ,
X_{\overline{h}})$, if we put $$ \pounds_{_{X^c}}dx^h=\alpha^h_m
dx^m + \beta^h_m \delta y^m.$$ Then we have
$$0=\pounds_{_{X^c}}(dx^h(X_i))=(\pounds_{_{X^c}}dx^h)X_i+dx^h(\pounds_{_{X^c}}X_i)=\alpha
^h_i-\partial_iv^h,$$ and $$
0=\pounds_{_{X^c}}(dx^h(X_{\overline{i}}))=(\pounds_{_{X^c}}dx^h)X_{
\overline{i}}+dx^h(\pounds_{_{X^c}}X_{\overline{i}})=\beta^h_i.$$Thus
we get (3). By the same way as the proof of (3), we can prove (4).
\end{proof}

\begin{lem}: Let $(M,g)$ be a Finsler manifold with Cartan
connection, then we have;\\
$(1)   \pounds_{_{X^c}}(g_{ij}dx^idx^j)=(\pounds_{_{V}}g_{ij})dx^i
dx^j ,$\newline $(2)   \pounds_{_{X^c}}(g_{ij}dx^i \delta
y^j)=g_{mi}(\pounds_{_{V}}N^m_{\,\,j})dx^idx^j+(\pounds_{_{V}}g_{ij})dx^i
\delta y^j,$\newline $(3)   \pounds_{_{X^c}}(g_{ij}dx^i \delta
y^j)=2(g_{mi}\pounds_{_{V}}N^m_{\,\,j})dx^i \delta
y^j+(\pounds_{_{V}}g_{ij})\delta y^i\delta y^j.$
\end{lem}
\begin{proof} By mean of above lemma, we get
\newline
$\pounds_{_{X^c}}(g_{ij}dx^idx^j)=X^c(g_{ij})dx^idx^j+2g_{ij}(\pounds_{_{X^c}}
dx^i)dx^j$ \newline
${\,\,\,\,\,\,\,\,\,\,\,\,\,\,\,\,\,\,\,\,\,\,\,\,\,\,\,\,\,\,\,\,\,\,\,\,\,\,\,\,}=
(v^hX_h+y^mv^h_{\,\,\,\,|m}X_{\overline{h}})
(g_{ij})dx^idx^j+2g_{ij}(\partial_mv^idx^m)dx^j$\newline
${\,\,\,\,\,\,\,\,\,\,\,\,\,\,\,\,\,\,\,\,\,\,\,\,\,\,\,\,\,\,\,\,\,\,\,\,\,\,\,\,}=
(\pounds_{_V}g_{ij})dx^idx^j.$\newline Thus we have (1). (2) and (3)
are easily proof by the same way as the proof of (1).
\end{proof}

 \setcounter{equation}{0}
 \setcounter{thm}{0}
\begin{thm} Let M be a Finsler manifold with Cartan connection and (TM,\~g)
its Riemannian (or Pseudo-Riemannian)  tangent bundle. Then every
conformal complete lift vector field on $TM$ is  homothetic.\\

\end{thm}
\begin{proof} Let V be a vector field on M, $X^c$ the complete lift vector field of V
which is conformal
 and \~g be a Pseudo-Riemannian
metric on TM derived from g. We have by definition $\pounds_{_{
X^c}}\tilde {g} =2\Omega \tilde {g} $. The Lie derivative of $\tilde
{g}$ gives
\begin{eqnarray}
\pounds_{_{ X^c}}\tilde g = \alpha (\pounds_{_{V }}g_{ij})dx^idx^j
+2\beta(\pounds_{_{V}} g_{ij})dx^i\delta y^j +2\beta g_{ai}
(\pounds_{_{V}}N^{a} _{j})dx^idx^j
\nonumber\\
+\gamma (\pounds_{_{V}}g_{ij})\delta y^i \delta y^j
 +2\gamma g_{aj}(\pounds_{_{V}}N^{a} _{i}) dx^i \delta y^j.
  \nonumber \\
   \end{eqnarray}
  So we have \begin{eqnarray}
  \pounds_{_{X^c}}\tilde g=[\alpha \pounds_{_V}g_{ij}  \nonumber
  +2\beta g_{ai} \pounds_{_V}N^{a} _{j} ]dx^idx^j \\ \nonumber
   +[2\beta\pounds_{_V} g_{ij} \nonumber
   +2\gamma g_{aj}\pounds_{_V}N^{a} _{i}  \nonumber
   ]dx^i\delta y^j\\ \nonumber
   +\gamma (\pounds_{_V}g_{ij})\delta y^i \delta y^j \\ \nonumber
  = 2\Omega {\tilde g}.\nonumber
\end{eqnarray}
Comparing with the definition of ${\tilde g}$, we find;
\begin{eqnarray}
 \alpha \pounds_{_{V}}g_{ij}+\beta(g_{ai}\pounds_{_V}N^{a} _{j}+
  g_{aj}\pounds_{_V}N^{a} _{i})=2\alpha\Omega g_{ij}.\\
 \beta\pounds_{_{V}}g_{ij}+\gamma g_{aj}\pounds_{_{V}}N^a_i=2\beta\Omega g_{ij}.    \\
\gamma \pounds_{_{V}}g_{ij} = 2\gamma\Omega g_{ij}.
 \end{eqnarray}
I) If $\gamma\neq0$ then from $(3.4)$ we have
$$\pounds_{_V}g _{ij}=2\Omega g_{ij},$$
and from (3.3) we have $$\pounds_{_{V}}N^a_i=0.$$ Using this and
$N_{i}^{h}=y^{m}F_{m \ i}^{h}$ we get
\begin{eqnarray}
 0=\pounds_{_{V}}
N_{i}^{h}=\pounds_{_{V}}(y^{m}F_{m \
i}^{\,h})=y^{m}\pounds_{_{V}}F_{\,m \ i}^{h}.
 \end{eqnarray}
 Where the last equality holds from equation (\ref{Lee der y}).\\
 II) If $\gamma=0$ since $\alpha\gamma-\beta^2\neq0$ we have
$\beta\neq0$ so from $(3.3)$ we have
$$\pounds_{_V}g _{ij}=2\Omega g_{ij},$$
and from (3.2) we have
$$g_{ai}\pounds_{_{V}} N_{j}^{a}+g_{aj}\pounds_{_{V}}N_{i}^{a}=0.$$
Using this and  equation (\ref{Lee der y}) and $N_{i}^{a}=y^{k}F_{k
\ i}^{a},$ we have
\begin{eqnarray}
y^k(g_{ai}\pounds_{_{V}} F_{k \ j}^{a}+ g_{aj}\pounds_{_{V}}F_{k \
i}^{a})=0.
\end{eqnarray}
In each case I) and II) we have
\begin{eqnarray}
 \pounds_{_V}g_{ij}=2\Omega g_{ij},
 \end{eqnarray} or from equation (1.6)
$$v^m\partial_mg_{ij}+g_{mj}\partial_iv^m+g_{im}\partial_jv^m+y^a\partial_av^m\dot{\partial
_m }g_{ij}=2\Omega g_{ij}.$$ Applying $\dot{\partial_k}$ to the both
side of the above equation, we find;
\begin{eqnarray} \nonumber
 2v^m\partial_mC_{ijk}+2C_{mjk}\partial_iv^m+2C_{imk}\partial_jv^m+
 2\partial_kv^mC_{ijm}+
 2y^a\partial_av^m\dot{\partial_k}C_{ijm}.\\ \nonumber=
 2g_{ij}\dot{\partial_k}\Omega+4\Omega C_{ijk}.
 \end{eqnarray}
 By using
$y^iC_{ijk}=0$, we obtain $\dot{\partial_k}\Omega =0$. Therefore
$\Omega$ is a function of $x$ alone.\newline
 From (\ref{interchange}) we have
  $$ y^k(\nabla_{k}\pounds_{_V} g_{ij}- \pounds_{_V}\nabla_{k}
 g_{ij})=y^k( g_{ai}\pounds_{_V}F^{a}_{jk}+g_{aj}\pounds_{_{V}}F^{a}
 _{ik}).$$
By using (3.5),(3.6) and (3.7) in each case I) and II) we find that
$$y^k\nabla_{k}\Omega=0.$$
Since $\Omega$ is a function of $x$ alone, we obtain $\partial_i
\Omega=0$. This together with connectedness of $M$, shows that
$\Omega$ is constant.
\end{proof}

Behroz Bidabad\\ Faculty of Mathematics, Amirkabir University of
Technology, (Tehran Polytechnic), 424, Hafez Ave. 15914,
  Tehran, Iran.

{\small\it E-mail: \textbf{bidabad@aut.ac.ir}}
\end{document}